\documentclass[a4paper,12pt]{article} 
\usepackage{amsmath} 
\usepackage{amsopn} 
\usepackage{amssymb} 
\usepackage[mathscr]{eucal} 
\usepackage{theorem} 
\usepackage{enumerate} 

\theorembodyfont{\itshape} 
\theoremstyle{plain}
  \newtheorem{thm}{Theorem}[section] 
  \newtheorem{pro}[thm]{Proposition}

\theorembodyfont{\rmfamily} 
\theoremstyle{plain}

\renewcommand{\theequation}%
           {\thesection.\arabic{equation}}

\begin{document} 

\begin{center} 
{\LARGE Time-like surfaces with zero mean curvature vector} 

\vspace{2mm} 

{\LARGE in 4-dimensional neutral space forms} 

\vspace{6mm} 

{\Large Naoya {\sc Ando}} 
\end{center} 

\vspace{3mm} 

\begin{quote} 
{\footnotesize \it Abstract} \ 
{\footnotesize 
Let $M$ be a Lorentz surface and 
$F:M\rightarrow N$ a time-like and conformal immersion of $M$ 
into a 4-dimensional neutral space form $N$ with zero mean curvature vector. 
We see that the curvature $K$ of the induced metric on $M$ by $F$ is 
identically equal to the constant sectional curvature $L_0$ of $N$ 
if and only if the covariant derivatives of 
both of the time-like twistor lifts are zero or light-like. 
If $K\equiv L_0$, then the normal connection $\nabla^{\perp}$ of $F$ is flat, 
while the converse is not necessarily true. 
We see that a holomorphic paracomplex quartic differential $Q$ on $M$ 
defined by $F$ is zero or null if and only if 
the covariant derivative of at least one of the time-like twistor lifts is 
zero or light-like. 
In addition, we see that $K$ is identically equal to $L_0$ 
if and only if not only $\nabla^{\perp}$ is flat 
but also $Q$ is zero or null.}
\end{quote} 

\vspace{3mm} 

\section{Introduction}\label{sect:intro} 

\setcounter{equation}{0} 

Let $N$ be a 4-dimensional neutral space form. 
Let $M$ be a Lorentz surface 
and $F:M\rightarrow N$ a time-like and conformal immersion of $M$ into $N$ 
with zero mean curvature vector. 
Then a holomorphic paracomplex quartic differential $Q$ on $M$ is 
defined by $F$. 
Isotropicity of $F$ is characterized by $Q\equiv 0$. 
If a time-like twistor lift of $F$ is horizontal, then $F$ is isotropic. 
However, although $F$ is isotropic, 
it is possible that the covariant derivatives of the time-like twistor lifts 
are not zero but light-like (\cite{ando4}). 
The covariant derivatives of 
both of the time-like twistor lifts are zero or light-like 
if and only if $F$ has one of the following properties: 
\begin{itemize} 
\item[{\rm (i)}]{the shape operator of a light-like normal vector field 
vanishes;} 
\item[{\rm (ii)}]{the shape operator of any normal vector field is 
zero or light-like} 
\end{itemize} 
(\cite{ando5}). 
See \cite{ando4}, \cite{ando5} for properties (i), (ii) of $F$ respectively. 
If $F$ satisfies (i) or (ii), 
then the curvature $K$ of the induced metric on $M$ by $F$ is 
identically equal to the constant sectional curvature $L_0$ of $N$. 
We will see that the converse holds, that is, 
if $K$ is identically equal to $L_0$, 
then the covariant derivatives of both of the time-like twistor lifts are 
zero or light-like (Theorem~\ref{thm:K=L_0}). 
In addition, we will study the following two properties related to $F$: 
\begin{itemize} 
\item[{\rm (I)}]{the normal connection $\nabla^{\perp}$ of $F$ is flat, 
that is, the curvature tensor $R^{\perp}$ of $\nabla^{\perp}$ vanishes;} 
\item[{\rm (II)}]{$Q$ is zero or null.} 
\end{itemize} 
We will see that if $K\equiv L_0$, then $\nabla^{\perp}$ is flat 
(Proposition~\ref{pro:R=0}), 
while the converse is not necessarily true. 
We will find time-like surfaces in $N$ with zero mean curvature vector, 
flat normal connection and $K\not= L_0$ (Theorem~\ref{thm:mu1vmu2u}). 
We will see that $Q$ is zero or null if and only if 
the covariant derivative of at least one of the time-like twistor lifts 
of $F$ is zero or light-like (Theorem~\ref{thm:Qzeronull}). 
In addition, we will find time-like surfaces in $N$ 
with zero mean curvature vector such that 
the covariant derivative of just one time-like twistor lift is 
zero or light-like (Theorem~\ref{thm:H=0Qzn}). 
We will see that $K$ is identically equal to $L_0$ 
if and only if $F$ satisfies both (I) and (II), that is, 
not only $\nabla^{\perp}$ is flat but also $Q$ is zero or null 
(Theorem~\ref{thm:K=L_0}). 

\vspace{3mm} 

\par\noindent 
\textit{Remark} \ 
Let $N$ be a 4-dimensional Riemannian, Lorentzian or neutral space form. 
Let $M$ be a Riemann or Lorentz surface. 
Let $F:M\rightarrow N$ be a space-like or time-like, and conformal immersion 
with zero mean curvature vector. 
Here $F$ is space-like or time-like, 
according to whether $M$ is Riemann or Lorentz. 
Let $K$ be the curvature of the induced metric by $F$ 
and $L_0$ the constant sectional curvature of $N$. 
If $N$ is either Riemannian or neutral, and if $M$ is Riemann, 
then by the equation of Gauss, 
a condition $K\equiv L_0$ means that $F$ is totally geodesic, and then 
by the equation of Ricci, the normal connection is flat. 
Suppose that $N$ is Lorentzian and that $M$ is Riemann. 
Then the condition $K\equiv L_0$ does not necessarily imply 
that $F$ is totally geodesic. 
By the equations of Gauss-Codazzi-Ricci, 
$K\equiv L_0$ is equivalent to a condition 
that the shape operator of a light-like normal vector field of $F$ vanishes, 
that is, 
a holomorphic quartic differential $Q$ on $M$ defined by $F$ vanishes 
(refer to \cite{ando3}). 
Then by the equation of Ricci, 
$K\equiv L_0$ means that the normal connection is flat. 
If $N$ is neutral and if $M$ is Lorentz, 
then we will obtain Proposition~\ref{pro:R=0} and Theorem~\ref{thm:K=L_0} 
below. 
If $N$ is Lorentzian and if $M$ is Lorentz, 
then $K\equiv L_0$ does not necessarily imply 
that $F$ is totally geodesic, 
while $K\equiv L_0$ means that the normal connection is flat 
(see the appendix below). 
We can refer to \cite{chen} 
for time-like surfaces in indefinite space forms 
with zero mean curvature vector and $K\equiv L_0$. 

\vspace{3mm} 

\par\noindent 
\textit{Remark} \ 
Let $N$ be an oriented $4$-dimensional Riemannian, Lorentzian or neutral 
manifold. 
Let $M$, $F$ be as in the beginning of the previous remark. 
Let $Q$ be a complex or paracomplex quartic differential on $M$ 
defined by $F$. 
If $N$ is Riemannian or neutral, and if $M$ is Riemann, 
then isotropicity of $F$ is equivalent to $Q\equiv 0$ 
and characterized by horizontality of a (space-like) twistor lift of $F$ 
(\cite{friedrich}, \cite{bryant}, \cite{ando4}). 
We can refer to \cite{ES}, \cite{BDM} 
for twistor spaces and space-like twistor spaces respectively. 
Horizontality of a (space-like) twistor lift is equivalent to 
parallelism of the complex structure corresponding to the lift, 
which is equivalent to a relation 
between the complex structure and the second fundamental form (\cite{ando7}). 
Suppose that $N$ is neutral and that $M$ is Lorentz. 
Then isotropicity of $F$ is equivalent to $Q\equiv 0$ 
and horizontality of a time-like twistor lift means $Q\equiv 0$, 
but $Q\equiv 0$ does not necessarily mean the horizontality (\cite{ando4}). 
We can refer to \cite{HM}, \cite{JR} for time-like twistor spaces. 
Horizontality of a time-like twistor lift is equivalent to 
parallelism of the paracomplex structure corresponding to the lift, 
which is equivalent to a relation 
between the paracomplex structure and the second fundamental form 
(\cite{ando7}). 
If $Q\equiv 0$ but if the time-like twistor lifts are not horizontal, 
then the covariant derivatives of them are zero or light-like (\cite{ando4}). 
Suppose that $N$ is Lorentzian. 
Based on the understanding of isotropicity in the above cases, 
the author introduced the definition of isotropicity of $F$ in \cite{ando6}. 
If $M$ is Riemann, 
then the isotropicity just means either $Q\equiv 0$ or 
a relation between a mixed-type structure given by a lift of $F$ and 
the second fundamental form with respect to 
a special local complex coordinate 
on a neighborhood of each point of $M$ (\cite{ando6}). 
In addition, 
if $N$ is a $4$-dimensional Lorentzian space form 
and if $F$ has no totally geodesic points, 
then the relation is equivalent to $Q\equiv 0$ (\cite{ando6}). 
If $M$ is Lorentz, 
then the isotropicity of $F$ is given by 
a relation between a mixed-type structure given by a lift of $F$ and 
the second fundamental form with respect to 
a special local paracomplex coordinate 
on a neighborhood of each point of $M$ (\cite{ando6}). 
In addition, 
if $N$ is a $4$-dimensional Lorentzian space form, 
then the isotropicity implies that $F$ is totally geodesic (\cite{ando6}), 
which is equivalent to $Q\equiv 0$. 

\vspace{3mm} 

\par\noindent 
\textit{Remark} \ 
Let $N$ be an oriented $4$-dimensional neutral manifold.  
Let $M$ be a Lorentz surface 
and $F:M\rightarrow N$ a time-like and conformal immersion 
with zero mean curvature vector. 
If the covariant derivative of a time-like twistor lift $\Theta$ of $F$ is 
zero or light-like, and not identically equal to zero at any point, 
then the covariant derivative of the paracomplex structure $J$ 
corresponding to $\Theta$ is locally represented as 
the tensor product of a nowhere zero $1$-form and a nilpotent structure 
related to $J$ (\cite{ando8}). 
We can refer to \cite{ando9}, \cite{ando8}, \cite{AK} 
for nilpotent structures. 
An oriented $4n$-dimensional neutral manifold equipped with 
such an almost paracomplex structure as $J$ becomes 
a Walker manifold (\cite{ando8}). 
We can find examples of such almost paracomplex structures 
on $E^{4n}_{2n}$ (\cite{ando8}). 
In addition, we can find all the pairs of sections of 
the two time-like twistor spaces associated with $E^4_2$ such that 
the covariant derivatives are zero or light-like, 
and not identically equal to zero at any point (\cite{ando8}). 

\section{Time-like surfaces with zero mean curvature vector}\label{sect:tsH=0} 

\setcounter{equation}{0} 

Let $N$ be a 4-dimensional neutral space form 
with constant sectional curvature $L_0$. 
If $L_0 =0$, 
then we can suppose $N=E^4_2$; 
if $L_0 >0$, 
then we can suppose $N=\{ x\in E^5_2 \ | \ \langle x, x\rangle =1/L_0 \}$; 
if $L_0 <0$, 
then we can suppose $N=\{ x\in E^5_3 \ | \ \langle x, x\rangle =1/L_0 \}$, 
where $\langle \ , \ \rangle$ is the metric of $E^5_2$ or $E^5_3$.  
Let $h$ be the neutral metric of $N$. 
Let $M$ be a Lorentz surface 
and $F:M\rightarrow N$ a time-like and conformal immersion of $M$ into $N$ 
with zero mean curvature vector. 
Let $\check{w} =u+jv$ be a local paracomplex coordinate of $M$, 
where we denote by $j$ the paraimaginary unit of paracomplex numbers. 
We denote by $\overline{\check{w}}$ the conjugate of $\check{w}$: 
$\overline{\check{w}} =u-jv$. 
We represent the induced metric $g$ on $M$ by $F$ 
as $g=e^{2\lambda} d\check{w}d\overline{\check{w}} 
     =e^{2\lambda} (du^2 -dv^2 )$ 
for a real-valued function $\lambda$. 
Then we have 
$$h(dF(\partial_{\check{w}} ), dF(\partial_{\overline{\check{w}}} )) 
 =\dfrac{e^{2\lambda}}{2} \quad 
  \left( \partial_{\check{w}} 
       :=\dfrac{1}{2} 
         \left( \dfrac{\partial}{\partial u} 
              +j\dfrac{\partial}{\partial v} 
         \right) , \ 
         \partial_{\overline{\check{w}}} 
       :=\dfrac{1}{2} 
         \left( \dfrac{\partial}{\partial u} 
              -j\dfrac{\partial}{\partial v} 
         \right) 
  \right) .$$ 
Let $\Tilde{\nabla}$ denote the connection of $E^4_2$, $E^5_2$ or $E^5_3$ 
according to $L_0 =0$, $>0$ or $<0$. 
Since $F$ has zero mean curvature vector, 
we have $\Tilde{\nabla}_{\partial /\partial \overline{\check{w}}} 
                      dF(\partial /\partial \check{w} ) 
       =-L_0 e^{2\lambda} F/2$. 
We set $T_1 :=dF(\partial /\partial u )$, $T_2 :=dF(\partial /\partial v )$. 
Let $N_1$, $N_2$ be normal vector fields of $F$ satisfying 
$$h(N_1 , N_1 )=-h(N_2 , N_2 )=e^{2\lambda} , \quad 
  h(N_1 , N_2 )=0.$$ 
Then we have 
\begin{equation} 
\begin{split} 
\Tilde{\nabla}_{T_1} (T_1 \ T_2 \ N_1 \ N_2 \ F) 
& =                  (T_1 \ T_2 \ N_1 \ N_2 \ F)S, \\ 
\Tilde{\nabla}_{T_2} (T_1 \ T_2 \ N_1 \ N_2 \ F) 
& =                  (T_1 \ T_2 \ N_1 \ N_2 \ F)T, 
\end{split} 
\label{dt1t2} 
\end{equation} 
where 
\begin{equation} 
S =\left[ \begin{array}{ccccc} 
       \lambda_u         & \lambda_v & -\alpha_1  &  \beta_1   & 1 \\ 
       \lambda_v         & \lambda_u &  \alpha_2  & -\beta_2   & 0 \\ 
       \alpha_1          & \alpha_2  &  \lambda_u &  \mu_1     & 0 \\ 
       \beta_1           & \beta_2   &  \mu_1     &  \lambda_u & 0 \\ 
       -L_0 e^{2\lambda} &  0        &   0        &   0        & 0 
            \end{array} 
   \right] , \quad 
T =\left[ \begin{array}{ccccc} 
       \lambda_v & \lambda_u         & -\alpha_2  &  \beta_2   & 0 \\ 
       \lambda_u & \lambda_v         &  \alpha_1  & -\beta_1   & 1 \\ 
       \alpha_2  & \alpha_1          &  \lambda_v &  \mu_2     & 0 \\ 
       \beta_2   & \beta_1           &  \mu_2     &  \lambda_v & 0 \\ 
        0        &  L_0 e^{2\lambda} &   0        &   0        & 0 
            \end{array} 
  \right] , 
\label{ST} 
\end{equation} 
and $\alpha_k$, $\beta_k$, $\mu_k$ ($k=1, 2$) are real-valued functions. 
From \eqref{dt1t2}, we obtain $S_v -T_u =ST-TS$. 
This is equivalent to the system of the equations of Gauss, Codazzi and Ricci. 
The equation of Gauss is given by 
\begin{equation} 
 \lambda_{uu} -\lambda_{vv} +L_0 e^{2\lambda} 
=\alpha^2_1 -\alpha^2_2 -(\beta^2_1 -\beta^2_2). 
\label{gauss}
\end{equation}
The equations of Codazzi are given by 
\begin{equation} 
\begin{split}
  (\alpha_1 )_v -(\alpha_2 )_u 
& =\alpha_2 \lambda_u -\alpha_1 \lambda_v 
  +\beta_2  \mu_1     -\beta_1  \mu_2  , \\ 
  (\alpha_2 )_v -(\alpha_1 )_u 
& =\alpha_1 \lambda_u -\alpha_2 \lambda_v 
  +\beta_1  \mu_1     -\beta_2  \mu_2  , \\ 
  (\beta_1 )_v -(\beta_2 )_u 
& =\beta_2 \lambda_u -\beta_1 \lambda_v 
  +\alpha_2  \mu_1     -\alpha_1  \mu_2  , \\ 
  (\beta_2 )_v -(\beta_1 )_u 
& =\beta_1 \lambda_u -\beta_2 \lambda_v 
  +\alpha_1  \mu_1     -\alpha_2  \mu_2 .  
\end{split} 
\label{codazzi} 
\end{equation} 
The equation of Ricci is given by 
\begin{equation} 
  (\mu_1 )_v -(\mu_2 )_u =2\alpha_1 \beta_2 -2\alpha_2 \beta_1 . 
\label{ricci}
\end{equation}
We set 
\begin{equation} 
X_{\pm} :=\alpha_1 \pm \beta_2 , \quad 
Y_{\pm} :=\alpha_2 \pm \beta_1 
\label{XpmYpm0} 
\end{equation} 
and 
\begin{equation} 
\phi_{\pm} :=\lambda_u \mp \mu_2  , \quad  
\psi_{\pm} :=\lambda_v \mp \mu_1 . 
\label{phps0} 
\end{equation} 
Then the system \eqref{codazzi} is rewritten into 
\begin{equation} 
\begin{split}
  (X_{\pm} )_v -(Y_{\pm} )_u & =\phi_{\pm} Y_{\pm} -\psi_{\pm} X_{\pm} , \\ 
  (Y_{\pm} )_v -(X_{\pm} )_u & =\phi_{\pm} X_{\pm} -\psi_{\pm} Y_{\pm} .  
\end{split}
\label{codazzi2}
\end{equation}
We set 
\begin{equation*} 
s:=\dfrac{1}{\sqrt{2}} (u+v), \quad 
t:=\dfrac{1}{\sqrt{2}} (u-v). 
\end{equation*} 
Then we have 
\begin{equation} 
 \dfrac{\partial}{\partial u} 
=\dfrac{1}{\sqrt{2}} 
 \left( \dfrac{\partial}{\partial s} +\dfrac{\partial}{\partial t} \right) , 
 \quad 
 \dfrac{\partial}{\partial v} 
=\dfrac{1}{\sqrt{2}} 
 \left( \dfrac{\partial}{\partial s} -\dfrac{\partial}{\partial t} \right) . 
\label{dudv} 
\end{equation} 
Therefore \eqref{codazzi2} is rewritten into 
\begin{equation} 
\begin{split}
(X_{\pm} +Y_{\pm} )_t 
& =-\dfrac{1}{\sqrt{2}} (\phi_{\pm} -\psi_{\pm})(X_{\pm} +Y_{\pm} ), \\ 
(X_{\pm} -Y_{\pm} )_s 
& =-\dfrac{1}{\sqrt{2}} (\phi_{\pm} +\psi_{\pm})(X_{\pm} -Y_{\pm} ). 
\end{split}
\label{codazzi3}
\end{equation}
A system which consists of \eqref{gauss} and \eqref{ricci} is 
rewritten into 
\begin{equation} 
\lambda_{uu} -\lambda_{vv} +L_0 e^{2\lambda} 
\pm ((\mu_1 )_v -(\mu_2 )_u ) 
=X^2_{\pm} -Y^2_{\pm} . 
\label{gaussricci} 
\end{equation} 
Applying \eqref{phps0} to \eqref{gaussricci}, we obtain 
\begin{equation}
 X^2_{\pm} -Y^2_{\pm} 
=(\phi_{\pm} )_u -(\psi_{\pm} )_v +L_0 e^{2\lambda} . 
\label{ricci3}
\end{equation}

Suppose that $F$ satisfies (i) in Section~\ref{sect:intro}, that is, 
the shape operator of a light-like normal vector field vanishes. 
Then we have $\beta_k =\varepsilon \alpha_k$ 
for $k=1, 2$ and $\varepsilon \in \{ +, -\}$. 
From \eqref{gauss}, we have 
\begin{equation} 
\lambda_{uu} -\lambda_{vv} +L_0 e^{2\lambda} =0. 
\label{K=L0} 
\end{equation} 
From \eqref{ricci}, we have $\gamma_u =\mu_1$, $\gamma_v =\mu_2$ 
for a function $\gamma$ of two variables $u$, $v$. 
From \eqref{codazzi}, we have 
\begin{equation*} 
(\alpha_+ )_t +\alpha_+ \lambda_t +\varepsilon \alpha_+ \gamma_t =0, \quad 
(\alpha_- )_s +\alpha_- \lambda_s +\varepsilon \alpha_- \gamma_s =0, 
\label{codazzi6} 
\end{equation*} 
where $\alpha_{\pm} :=\alpha_1 \pm \alpha_2$. 
Then there exist functions $p_{\pm}$ of one variable 
satisfying $\alpha_{\pm} =p_{\pm} (u\pm v)e^{-\lambda -\varepsilon \gamma}$. 
Therefore we obtain 
\begin{equation*} 
\begin{split} 
\alpha_1 & = 
\dfrac{1}{2} (p_+ (u+v)+p_- (u-v))e^{-\lambda -\varepsilon \gamma} , \\ 
\alpha_2 & = 
\dfrac{1}{2} (p_+ (u+v)-p_- (u-v))e^{-\lambda -\varepsilon \gamma} . 
\end{split} 
\label{alpha} 
\end{equation*} 

Suppose that $F$ satisfies (ii) in Section~\ref{sect:intro}, that is, 
the shape operator of any normal vector field is zero or light-like. 
Then, as in \cite{ando5}, 
we obtain \eqref{K=L0} and 
$$\alpha_2 =\varepsilon \alpha_1 , \quad 
  \beta_2  =\varepsilon \beta_1  , \quad 
  \mu_1    =            \gamma_u , \quad 
  \mu_2    =            \gamma_v$$ 
for $\varepsilon \in \{ +, -\}$ and 
a function $\gamma$ of two variables $u$, $v$, 
and $\alpha_1$, $\beta_1$ are represented as 
\begin{equation*} 
\begin{split} 
\alpha_1 & = \dfrac{1}{2e^{\lambda}} (\phi (u+\varepsilon v)e^{\gamma} 
                                     +\psi (u+\varepsilon v)e^{-\gamma} ), \\ 
\beta_1  & =-\dfrac{1}{2e^{\lambda}} (\phi (u+\varepsilon v)e^{\gamma} 
                                     -\psi (u+\varepsilon v)e^{-\gamma} ) 
\end{split} 
\label{ab} 
\end{equation*} 
for functions $\phi$, $\psi$ of one variable. 

\section{Time-like twistor lifts}\label{sect:ttl} 

\setcounter{equation}{0} 

Let $N$, $h$, $M$, $F$ be as in the beginning of the previous section. 
Then the $2$-fold exterior power $\bigwedge^2\!F^*\!T\!N$ 
of the pull-back bundle $F^*\!T\!N$ on $M$ by $F$ is 
a vector bundle on $M$ of rank $6$. 
Let $\hat{h}$ be the metric of $\bigwedge^2\!F^*\!T\!N$ induced by $h$. 
Then $\hat{h}$ has signature $(2,4)$. 
Noticing 
the double covering $SO_0(2,2)\longrightarrow SO_0 (1,2)\times SO_0(1,2)$, 
we see that $\bigwedge^2\!F^*\!T\!N$ is decomposed into 
two subbundles $\bigwedge^2_{\pm}\!F^*\!T\!N$ of rank $3$. 
These are orthogonal to each other with respect to $\hat{h}$ 
and the restriction of $\hat{h}$ on each of them has signature $(1,2)$. 
The time-like twistor spaces associated with $F^*\!T\!N$ are 
fiber bundles on $M$ defined by 
$$U_-\!\left.\left(\textstyle\bigwedge^2_{\pm}\!F^*\!T\!N\right) 
   :=\left\{ \theta \in \textstyle\bigwedge^2_{\pm}\!F^*\!T\!N \ \right| \ 
             \hat{h} (\theta , \theta )=-1\right\} .$$ 
Let $(e_1 , e_2 , e_3 , e_4 )$ be an ordered local pseudo-orthonormal 
frame field of $F^*\!T\!N$ giving the orientation of $N$ 
such that $e_1$, $e_2$ (respectively, $e_3$, $e_4$) are 
space-like (respectively, time-like). 
Set 
\begin{equation*} 
\begin{split} 
  \Theta_{\pm , 1} & 
:=\dfrac{1}{\sqrt{2}} (e_1 \wedge e_2 \pm e_3 \wedge e_4 ), \\ 
  \Theta_{\pm , 2} & 
:=\dfrac{1}{\sqrt{2}} (e_1 \wedge e_3 \pm e_4 \wedge e_2 ), \\ 
  \Theta_{\pm , 3} & 
:=\dfrac{1}{\sqrt{2}} (e_1 \wedge e_4 \pm e_2 \wedge e_3 ). 
\end{split} 
\end{equation*} 
Then $\bigwedge^2_{\pm}\!F^*\!T\!N$ are locally generated 
by $\Theta_{\mp , 1}$, $\Theta_{\pm , 2}$, $\Theta_{\pm , 3}$, 
respectively. 
The time-like twistor lifts $\Theta_{F, \pm}$ of $F$ are sections 
of $U_-\!\left(\textstyle\bigwedge^2_{\pm}\!F^*\!T\!N\right)$ 
and locally represented as $\Theta_{F, \pm} =\Theta_{\pm , 2}$, 
where we suppose that $e_1$, $e_3$ satisfy $e_1$, $e_3 \in dF(T\!M)$ 
so that $(e_1 , e_3 )$ gives the orientation of $M$. 

Let $\nabla$ be the Levi-Civita connection of $N$. 
Then $\nabla$ induces a connection $\hat{\nabla}$ of $\bigwedge^2\!F^*\!T\!N$. 
In addition, $\hat{\nabla}$ gives connections 
of $\bigwedge^2_{\pm}\!F^*\!T\!N$. 
If we set 
\begin{equation*} 
e_1 :=\dfrac{1}{e^{\lambda}} T_1 , \quad 
e_3 :=\dfrac{1}{e^{\lambda}} T_2 , \quad 
e_2 :=\dfrac{1}{e^{\lambda}} N_1 , \quad 
e_4 :=\dfrac{1}{e^{\lambda}} N_2 , 
\end{equation*} 
then the covariant derivatives of $\Theta_{F, \pm} =\Theta_{\pm , 2}$ 
with respect to $\hat{\nabla}$ are given by 
\begin{equation*} 
\begin{split} 
   \hat{\nabla}_{T_1} \Theta_{F, \pm} 
& =Y_{\pm} \Theta_{\mp , 1} \pm X_{\pm} \Theta_{\pm , 3} , \\ 
   \hat{\nabla}_{T_2} \Theta_{F, \pm} 
& =X_{\pm} \Theta_{\mp , 1} \pm Y_{\pm} \Theta_{\pm , 3} , 
\end{split} 
\end{equation*} 
where $X_{\pm}$, $Y_{\pm}$ are as in \eqref{XpmYpm0}. 
Therefore we obtain 

\begin{pro}\label{pro:hatnabla} 
The following hold\/$:$ 
\begin{itemize} 
\item[{\rm (a)}]{one of $\hat{\nabla} \Theta_{F, \pm}$ is zero or light-like 
if and only if $F$ satisfies either $X^2_+ =Y^2_+$ or $X^2_- =Y^2_- ;$} 
\item[{\rm (b)}]{both of $\hat{\nabla} \Theta_{F, \pm}$ are zero or light-like 
if and only if $F$ satisfies both of $X^2_{\pm} =Y^2_{\pm}$.} 
\end{itemize} 
\end{pro} 

We see that $F$ satisfies $X^2_+ =Y^2_+$ 
if and only if $F$ satisfies one of 
\begin{equation} 
\alpha_1 +\alpha_2 =-(\beta_1 +\beta_2 ), \quad 
\alpha_1 -\alpha_2 =  \beta_1 -\beta_2 
\label{X^2_+=Y^2_+} 
\end{equation} 
and that $F$ satisfies $X^2_- =Y^2_-$ 
if and only if $F$ satisfies one of 
\begin{equation} 
\alpha_1 +\alpha_2 =  \beta_1 +\beta_2 , \quad 
\alpha_1 -\alpha_2 =-(\beta_1 -\beta_2 ). 
\label{X^2_-=Y^2_-} 
\end{equation} 
We set $\gamma_k :=\alpha_k +j\beta_k$ for $k=1, 2$. 
Then $F$ satisfies one relation in \eqref{X^2_+=Y^2_+} and \eqref{X^2_-=Y^2_-} 
if and only if one of $\gamma_1 \pm \gamma_2$ is zero or null. 
Let $\sigma$ be the second fundamental form of $F$. 
Set $\sigma_{\check{w} \check{w}} 
   :=\sigma (\partial /\partial \check{w} , 
             \partial /\partial \check{w} )$ 
and $Q:=h(\sigma_{\check{w} \check{w}} , 
          \sigma_{\check{w} \check{w}} )d\check{w}^4$. 
Then $Q$ is a paracomplex quartic differential defined on $M$. 
Since $N$ is a neutral space form, 
we see by the equations of Codazzi that $Q$ is paraholomorphic. 
From \eqref{dt1t2} with \eqref{ST}, we have 
\begin{equation*} 
 \sigma_{\check{w} \check{w}}  
=\dfrac{1}{2} (\alpha_1 N_1 +\beta_1 N_2 +j(\alpha_2 N_1 +\beta_2 N_2 )) 
\end{equation*} 
and therefore we obtain 
\begin{equation} 
Q=\dfrac{e^{2\lambda}}{4} (|\gamma_1 |^2 +|\gamma_2 |^2 
                        +2j{\rm Re}\,(\gamma_1 \overline{\gamma}_2 ))
 d\check{w}^4 , 
\label{Q} 
\end{equation} 
where $|\gamma_k |^2$ is the square norm of $\gamma_k$: 
$|\gamma_k |^2 :=\alpha^2_k -\beta^2_k$. 
From \eqref{Q} and 
\begin{equation*} 
 |\gamma_1 \pm \gamma_2 |^2 
=|\gamma_1 |^2 +|\gamma_2 |^2 \pm 2{\rm Re}\,(\gamma_1 \overline{\gamma}_2 ), 
\label{gamma1pmgamma2} 
\end{equation*} 
we see that one of $\gamma_1 \pm \gamma_2$ is zero or null 
if and only if $Q$ is zero or null. 
Therefore from (a) of Proposition~\ref{pro:hatnabla}, we obtain 

\begin{thm}\label{thm:Qzeronull} 
One of $\hat{\nabla} \Theta_{F, \pm}$ is zero or light-like 
if and only if $Q$ is zero or null. 
\end{thm} 

\section{Time-like surfaces with zero mean curvature vector 
and \mbox{\boldmath{$K\equiv L_0$}}}\label{sect:K=L0} 

\setcounter{equation}{0} 

Suppose that the curvature $K$ of the induced metric $g$ is 
identically equal to $L_0$. 
Then we have \eqref{K=L0}. 
By this and \eqref{gaussricci}, we obtain 
\begin{equation}
  (\mu_1 )_v -(\mu_2 )_u =\pm (X^2_{\pm} -Y^2_{\pm} ). 
\label{ricci2}
\end{equation}
From \eqref{ricci2}, we obtain 
\begin{equation}
X^2_+ -Y^2_+ =-(X^2_- -Y^2_- ). 
\label{XpmYpm}
\end{equation}
We will prove 

\begin{pro}\label{pro:R=0}
If $K\equiv L_0$, 
then $(\mu_1 )_v \equiv (\mu_2 )_u$. 
\end{pro}

\vspace{3mm}

\par\noindent 
\textit{Proof} \ 
Suppose $(\mu_1 )_v \not= (\mu_2 )_u$ on a neighborhood of a point of $M$. 
Then from \eqref{ricci2}, 
we have $X^2_+ \not= Y^2_+$ and $X^2_- \not= Y^2_-$. 
We set 
\begin{equation} 
P_{\pm} :=\log |X_{\pm} +Y_{\pm} |, \quad 
Q_{\pm} :=\log |X_{\pm} -Y_{\pm} |. 
\label{PQ} 
\end{equation} 
Then noticing \eqref{XpmYpm}, we obtain 
\begin{equation} 
P_+ +Q_+ =P_- +Q_- =\log |X^2_{\pm} -Y^2_{\pm} |. 
\label{PpmQpm} 
\end{equation} 
From \eqref{XpmYpm} and \eqref{PpmQpm}, we obtain 
\begin{equation} 
X^2_{\pm} -Y^2_{\pm} =\pm \varepsilon e^{P_{\pm} +Q_{\pm}} , 
\label{epsilon} 
\end{equation} 
where $\varepsilon \in \{ +, -\}$. 
It follows from \eqref{PQ} that \eqref{codazzi3} is rewritten into 
\begin{equation*} 
(P_{\pm} )_t =-\dfrac{1}{\sqrt{2}} (\phi_{\pm} -\psi_{\pm}), \quad 
(Q_{\pm} )_s =-\dfrac{1}{\sqrt{2}} (\phi_{\pm} +\psi_{\pm}). 
\label{codazzi4}
\end{equation*}
Therefore we obtain 
\begin{equation} 
\phi_{\pm} =-\dfrac{1}{\sqrt{2}} ((Q_{\pm} )_s +(P_{\pm} )_t ), \quad 
\psi_{\pm} =-\dfrac{1}{\sqrt{2}} ((Q_{\pm} )_s -(P_{\pm} )_t ). 
\label{codazzi5}
\end{equation} 
From \eqref{dudv} and \eqref{codazzi5}, we obtain 
\begin{equation}
(\phi_{\pm} )_u -(\psi_{\pm} )_v 
=-(P_{\pm} +Q_{\pm} )_{st} . 
\label{phu-psv}
\end{equation}
Applying \eqref{epsilon}, \eqref{phu-psv} to \eqref{ricci3}, we obtain 
\begin{equation}
(P_{\pm} +Q_{\pm} )_{st} \pm \varepsilon e^{P_{\pm} +Q_{\pm}} 
=L_0 e^{2\lambda} . 
\label{R}
\end{equation} 
However, \eqref{R} contradicts to \eqref{PpmQpm}. 
Therefore we can not suppose $(\mu_1 )_v \not= (\mu_2 )_u$. 
Hence we obtain $(\mu_1 )_v \equiv (\mu_2 )_u$. 
\hfill 
$\square$ 

\vspace{3mm} 

\par\noindent 
\textit{Remark} \ 
We see that $F$ satisfies $(\mu_1 )_v \equiv (\mu_2 )_u$ 
if and only if the curvature tensor $R^{\perp}$ of 
the normal connection $\nabla^{\perp}$ of $F$ vanishes. 
Although $F$ satisfies $(\mu_1 )_v \equiv (\mu_2 )_u$, 
$K$ is not necessarily equal to $L_0$ (see Section~\ref{sect:fnc}). 

\vspace{3mm} 

We will prove 

\begin{thm}\label{thm:K=L_0} 
Let $N$ be a $4$-dimensional neutral space form 
with constant sectional curvature $L_0$. 
Let $M$ be a Lorentz surface and 
$F:M\rightarrow N$ a time-like and conformal immersion of $M$ into $N$ 
with zero mean curvature vector. 
Then the following are mutually equivalent\/$:$ 
\begin{itemize}
\item[{\rm (a)}]{the curvature $K$ of the induced metric by $F$ is 
identically equal to $L_0 ;$} 
\item[{\rm (b)}]{the covariant derivatives of 
both of the time-like twistor lifts are zero or light-like\/$;$} 
\item[{\rm (c)}]{not only $\nabla^{\perp}$ is flat 
but also $Q$ is zero or null.} 
\end{itemize}
\end{thm}

\vspace{3mm} 

\par\noindent 
\textit{Proof} \ 
Suppose $K\equiv L_0$. 
Then by \eqref{gaussricci} and Proposition~\ref{pro:R=0}, 
we have both $X^2_+ =Y^2_+$ and $X^2_- =Y^2_-$. 
Therefore we see from (b) of Proposition~\ref{pro:hatnabla} that 
the covariant derivatives of both of the time-like twistor lifts of $F$ are 
zero or light-like. 
If the covariant derivatives of both of the time-like twistor lifts are 
zero or light-like, 
then from \eqref{gaussricci} and Proposition~\ref{pro:hatnabla}, 
we obtain $K\equiv L_0$. 
Hence we have proved that (a) and (b) are equivalent. 
By \eqref{gaussricci}, Proposition~\ref{pro:hatnabla}, 
Theorem~\ref{thm:Qzeronull} and Proposition~\ref{pro:R=0}, 
we obtain (c) from (a). 
In addition, by \eqref{gaussricci}, Proposition~\ref{pro:hatnabla} 
and Theorem~\ref{thm:Qzeronull}, we obtain (a) from (c).  
Hence we have proved that (a) and (c) are equivalent. 
\hfill 
$\square$ 

\vspace{3mm} 

\par\noindent 
\textit{Remark} \ 
Let $\lambda$ be a function of two variables $u$, $v$ 
satisfying \eqref{K=L0}. 
If $L_0 =0$, 
then $\lambda$ is represented as the sum of a function of $s$ 
and a function of $t$. 
If $L_0 \not= 0$, 
then \eqref{K=L0} is easily transformed to Liouville's equation: 
if $L_0 >0$, 
then $\lambda$ is represented as in the form of 
$$e^{2\lambda} 
 =\dfrac{2p'(\sqrt{L_0} s)q'(\sqrt{L_0} t)}{
        (p(\sqrt{L_0} s) -q(\sqrt{L_0} t))^2} ,$$ 
where $p$, $q$ are functions of one variable; 
if $L_0 <0$, 
then $\lambda$ is represented as in the form of 
$$e^{2\lambda} 
 =\dfrac{2p'(\sqrt{|L_0 |} s)q'(\sqrt{|L_0 |} t)}{
        (p(\sqrt{|L_0 |} s) +q(\sqrt{|L_0 |} t))^2} ,$$ 
where $p$, $q$ are as above. 

\section{Time-like surfaces with zero mean curvature vector 
and flat normal connection}\label{sect:fnc} 

\setcounter{equation}{0} 

Let $F$ be as in the beginning of Section~\ref{sect:tsH=0}. 
Suppose that the normal connection of $F$ is flat. 
Then $F$ satisfies $(\mu_1 )_v \equiv (\mu_2 )_u$. 
From \eqref{gaussricci}, we have 
\begin{equation} 
\lambda_{uu} -\lambda_{vv} +L_0 e^{2\lambda} 
=X^2_{\pm} -Y^2_{\pm} . 
\label{gaussricci2} 
\end{equation} 
From \eqref{gaussricci2}, we obtain 
\begin{equation}
X^2_+ -Y^2_+ =X^2_- -Y^2_- . 
\label{XpmYpm2}
\end{equation} 
Suppose $X^2_{\pm} \not= Y^2_{\pm}$, i.e., $K\not= L_0$. 
Let $P_{\pm}$, $Q_{\pm}$ be as in \eqref{PQ}. 
Then we have \eqref{PpmQpm}. 
From \eqref{PpmQpm} and \eqref{XpmYpm2}, we obtain 
\begin{equation} 
X^2_{\pm} -Y^2_{\pm} =\varepsilon e^R , 
\label{epsilon2} 
\end{equation} 
where $\varepsilon \in \{ +, -\}$ and $R:=P_{\pm} +Q_{\pm}$. 
By \eqref{codazzi3} and \eqref{PQ}, 
we obtain \eqref{codazzi5} and \eqref{phu-psv}. 
By \eqref{phps0} and \eqref{codazzi5}, we obtain 
\begin{equation} 
\begin{split} 
\lambda_u & =-\dfrac{1}{2\sqrt{2}} ((Q_+ +Q_- )_s +(P_+ +P_- )_t ), \\ 
\lambda_v & =-\dfrac{1}{2\sqrt{2}} ((Q_+ +Q_- )_s -(P_+ +P_- )_t ). 
\end{split} 
\label{lulv} 
\end{equation} 
Applying \eqref{lulv} to $\lambda_{uv} =\lambda_{vu}$ 
and using \eqref{dudv}, 
we obtain $(P_+ +P_- )_{st} =(Q_+ +Q_- )_{st}$. 
By this and \eqref{PpmQpm}, 
we see that there exist functions $x_1$, $y_1$ of one variable $s$, $t$ 
respectively satisfying 
\begin{equation} 
Q_- =P_+ +x_1 (s) +y_1 (t), \quad 
Q_+ =P_- +x_1 (s) +y_1 (t). 
\label{QP} 
\end{equation} 
Applying $R=P_{\pm} +Q_{\pm}$ to \eqref{phu-psv}, we obtain 
\begin{equation}
(\phi_{\pm} )_u -(\psi_{\pm} )_v 
=-R_{st} . 
\label{phu-psv2}
\end{equation}
By $(\mu_1 )_v \equiv (\mu_2 )_u$ and \eqref{phps0}, we have 
\begin{equation} 
(\phi_{\pm} )_u -(\psi_{\pm} )_v 
=\lambda_{uu} -\lambda_{vv} . 
\label{pupvl} 
\end{equation} 
From \eqref{gaussricci2}, \eqref{epsilon2}, \eqref{phu-psv2} 
and \eqref{pupvl}, we obtain 
\begin{equation} 
R_{st} =-2\lambda_{st} 
 =L_0 e^{2\lambda} -\varepsilon e^R . 
\label{R2} 
\end{equation} 
From the first equation in \eqref{R2}, 
we have 
\begin{equation} 
R=-2\lambda +x_2 (s) +y_2 (t) 
\label{x2sy2t} 
\end{equation} 
for functions $x_2$, $y_2$ of one variable. 
Applying this to the second equation in \eqref{R2}, we obtain 
\begin{equation} 
\lambda_{st} =-\dfrac{L_0}{2} e^{2\lambda} 
              +\dfrac{\varepsilon}{2} e^{x_2 (s)} e^{y_2 (t)} e^{-2\lambda} . 
\label{lambda} 
\end{equation} 
From $R=P_{\pm} +Q_{\pm}$ and \eqref{QP}, 
we obtain 
\begin{equation} 
P_+ +P_- =R-x_1 (s) -y_1 (t), \quad 
Q_+ +Q_- =R+x_1 (s) +y_1 (t). 
\label{P+P-Q+Q-} 
\end{equation} 
Applying \eqref{x2sy2t} and \eqref{P+P-Q+Q-} 
to \eqref{lulv}, we obtain 
\begin{equation*} 
x_2 (s)=-x_1 (s)+a, \quad 
y_2 (t)= y_1 (t)+b 
\label{ab2} 
\end{equation*} 
for constants $a$, $b$. 
Let $\Tilde{s}$, $\Tilde{t}$ be functions of $s$, $t$ respectively 
satisfying 
$$\dfrac{d\Tilde{s}}{ds} =e^{\frac{1}{2} x_2 (s)} , \quad 
  \dfrac{d\Tilde{t}}{dt} =e^{\frac{1}{2} y_2 (t)} .$$ 
Then $(\Tilde{s} , \Tilde{t} )$ are local coordinates 
and from \eqref{lambda}, we obtain 
\begin{equation*} 
  \Tilde{\lambda}_{\Tilde{s} \Tilde{t}} 
=-\dfrac{L_0}{2} e^{2\Tilde{\lambda}} 
 +\dfrac{\varepsilon}{2} e^{-2\Tilde{\lambda}} , 
\label{Tildelambda} 
\end{equation*} 
where $\Tilde{\lambda} :=\lambda -(1/4)(x_2 (s)+y_2 (t))$. 
We set $\Tilde{u} :=(1/\sqrt{2})(\Tilde{s} +\Tilde{t} )$, 
       $\Tilde{v} :=(1/\sqrt{2})(\Tilde{s} -\Tilde{t} )$. 
Then $(\Tilde{u} , \Tilde{v} )$ are local coordinates 
satisfying $g=e^{2\Tilde{\lambda}} (d\Tilde{u}^2 -d\Tilde{v}^2 )$. 
Let $\Tilde{X}_{\pm}$, $\Tilde{Y}_{\pm}$ be defined as in \eqref{XpmYpm0} 
for $(\Tilde{u} , \Tilde{v} )$. 
Then we have 
\begin{equation} 
 \Tilde{\lambda}_{\Tilde{u} \Tilde{u}} 
-\Tilde{\lambda}_{\Tilde{v} \Tilde{v}} +L_0 e^{2\Tilde{\lambda}} 
=\Tilde{X}^2_{\pm} -\Tilde{Y}^2_{\pm} . 
\label{gaussricci2Tilde} 
\end{equation} 
Comparing \eqref{gaussricci2} with \eqref{gaussricci2Tilde}, 
we obtain 
\begin{equation*} 
 \Tilde{X}^2_{\pm} -\Tilde{Y}^2_{\pm} 
=(X^2_{\pm} -Y^2_{\pm} )e^{-(1/2)(x_2 (s)+y_2 (t))} . 
\label{TildeXX} 
\end{equation*} 
Let $\Tilde{P}_{\pm}$, $\Tilde{Q}_{\pm}$ be defined as in \eqref{PQ} 
for $(\Tilde{u} , \Tilde{v} )$ 
and set $\Tilde{R} :=\Tilde{P}_{\pm} +\Tilde{Q}_{\pm}$. 
Then we have $\Tilde{R} =R-(1/2)(x_2 (s)+y_2 (t))$. 
Therefore by this and \eqref{x2sy2t}, 
we obtain $\Tilde{R} =-2\Tilde{\lambda}$. 

Let $\lambda$ be a function of two variables $s$, $t$ 
satisfying 
\begin{equation} 
\lambda_{st} =-\dfrac{L_0}{2} e^{2\lambda} 
              +\dfrac{\varepsilon}{2} e^{-2\lambda} 
\label{lambdast} 
\end{equation} 
for $\varepsilon \in \{ +, -\}$. 
We set $R:=-2\lambda$. 
Then $\lambda$, $R$ satisfy \eqref{R2}. 
Let $P_{\pm}$, $Q_{\pm}$ be functions satisfying 
\begin{equation} 
R=P_{\pm} +Q_{\pm} , \quad 
Q_- =P_+ +c, \quad 
Q_+ =P_- +c 
\label{PpmQpmc} 
\end{equation} 
for a constant $c$. 
Then $P_{\pm}$, $Q_{\pm}$ satisfy \eqref{lulv}. 
Let $\phi_{\pm}$, $\psi_{\pm}$ be functions as in \eqref{codazzi5}. 
Then from \eqref{lulv}, we obtain 
\begin{equation} 
\lambda_u -\phi_+ =-(\lambda_u -\phi_- ), \quad 
\lambda_v -\psi_+ =-(\lambda_v -\psi_- ). 
\label{mu0} 
\end{equation} 
Noticing \eqref{mu0}, 
we see that 
\begin{equation} 
\begin{split} 
\mu_1 & := (\lambda_v -\psi_+ =) \ 
            \dfrac{1}{2\sqrt{2}} ((Q_+ -Q_- )_s -(P_+ -P_- )_t ), \\ 
\mu_2 & := (\lambda_u -\phi_+ =) \ 
            \dfrac{1}{2\sqrt{2}} ((Q_+ -Q_- )_s +(P_+ -P_- )_t ) 
\end{split} 
\label{mu1mu2} 
\end{equation} 
satisfy \eqref{phps0}. 
In addition, noticing 
$$(\phi_+ )_u -(\psi_+ )_v =-R_{st} =\lambda_{uu} -\lambda_{vv} ,$$ 
we obtain $(\mu_1 )_v \equiv (\mu_2 )_u$. 
We set 
\begin{equation} 
X_{\pm} 
 :=\dfrac{\varepsilon'_{\pm}}{2} (e^{P_{\pm}} +\varepsilon e^{Q_{\pm}} ), 
   \quad 
Y_{\pm} 
 :=\dfrac{\varepsilon'_{\pm}}{2} (e^{P_{\pm}} -\varepsilon e^{Q_{\pm}} ), 
\label{setXpmYpm} 
\end{equation} 
where $\varepsilon'_+$, $\varepsilon'_- \in \{ +, -\}$. 
Then we have \eqref{epsilon2}. 
Therefore we see from \eqref{lambdast} 
that $\lambda$, $X_{\pm}$, $Y_{\pm}$ satisfy \eqref{gaussricci2}. 
We set 
\begin{equation} 
\begin{array}{lcl} 
\alpha_1 :=\dfrac{1}{2} (X_+ +X_- ), & \ & 
\alpha_2 :=\dfrac{1}{2} (Y_+ +Y_- ), \\ 
\ & \ & \\ 
\beta_1  :=\dfrac{1}{2} (Y_+ -Y_- ), & \ & 
\beta_2  :=\dfrac{1}{2} (X_+ -X_- ).   
  \end{array} 
\label{alphabeta}
\end{equation} 
Then $S$, $T$ with \eqref{ST} satisfy $S_v -T_u =ST-TS$. 

Hence we obtain 

\begin{thm}\label{thm:mu1vmu2u} 
Let $N$ be a 4-dimensional neutral space form 
with constant sectional curvature $L_0$. 
Let $M$ be a Lorentz surface 
and $F:M\rightarrow N$ a time-like and conformal immersion of $M$ into $N$ 
with zero mean curvature vector.  
Then the following are equivalent\/$:$ 
\begin{itemize}
\item[{\rm (a)}]{$F$ satisfies both $K\not= L_0$ and $R^{\perp} =0;$} 
\item[{\rm (b)}]{the induced metric $g$ by $F$ is locally represented 
as $g=e^{2\lambda} d\check{w}d\overline{\check{w}}$ 
for a real-valued function $\lambda$ satisfying 
\begin{equation} 
 \lambda_{uu} -\lambda_{vv} 
=-L_0 e^{2\lambda} +\varepsilon e^{-2\lambda} \quad 
(\varepsilon \in \{ +, -\} ) 
\label{lambdaeq} 
\end{equation} 
and $\alpha_k$, $\beta_k$, $\mu_k$ $(k=1, 2)$ satisfy 
\eqref{mu1mu2}, \eqref{alphabeta} 
with \eqref{PpmQpmc}, \eqref{setXpmYpm} and $R=-2\lambda$.} 
\end{itemize} 
In addition, a function $\lambda$ of two variables $u$, $v$ 
satisfying \eqref{lambdaeq} 
and functions $\alpha_k$, $\beta_k$, $\mu_k$ $(k=1, 2)$ of $u$, $v$ 
given by \eqref{mu1mu2}, \eqref{alphabeta} 
with \eqref{PpmQpmc}, \eqref{setXpmYpm} and $R=-2\lambda$ 
locally define a time-like surface in $N$ 
with zero mean curvature vector, $K\not= L_0$ and $R^{\perp} =0$, 
which is unique up to an isometry of $N$. 
\end{thm} 

\section{Time-like surfaces with zero mean curvature vector 
and \mbox{\boldmath{$X^2_+ \not= Y^2_+$}}, 
    \mbox{\boldmath{$X_- =Y_- \not= 0$}}}\label{sect:H=0Qzn} 

\setcounter{equation}{0} 

Let $F$ be as in the beginning of Section~\ref{sect:tsH=0}. 
Suppose that the covariant derivative of 
just one of the time-like twistor lifts of $F$ is zero or light-like. 
In this section, 
we suppose that $F$ satisfies $X^2_+ \not= Y^2_+$ and $X_- =Y_- \not= 0$. 
Let $P_{\pm}$, $Q_+$ be as in \eqref{PQ}. 
Then $\phi_+$, $\psi_+$ are given by \eqref{codazzi5}, 
and therefore by \eqref{phps0}, $\phi_-$, $\psi_-$ are given by 
\begin{equation} 
\phi_- =2\lambda_u +\dfrac{1}{\sqrt{2}} ((Q_+ )_s +(P_+ )_t ), \quad 
\psi_- =2\lambda_v +\dfrac{1}{\sqrt{2}} ((Q_+ )_s -(P_+ )_t ). 
\label{phi-psi-} 
\end{equation}  
Noticing $P_- =\log 2|X_- |$, 
from one relation in \eqref{codazzi3} and \eqref{phi-psi-}, we obtain 
\begin{equation*} 
(P_- )_t =-2\lambda_t -(P_+ )_t . 
\label{PlambdaP} 
\end{equation*} 
Therefore we obtain 
\begin{equation} 
\lambda =-\dfrac{1}{2} (P_+ +P_- )+a(s), 
\label{lambdaPPa} 
\end{equation} 
where $a$ is a function of one variable. 
Applying \eqref{lambdaPPa} to \eqref{phi-psi-}, we obtain 
\begin{equation} 
\begin{split} 
\phi_-  & =-\dfrac{1}{\sqrt{2}} (P_+ -Q_+ )_s -(P_- )_u +\sqrt{2} a'(s), \\ 
\psi_-  & =-\dfrac{1}{\sqrt{2}} (P_+ -Q_+ )_s -(P_- )_v +\sqrt{2} a'(s). 
\end{split} 
\label{phi-psi-2} 
\end{equation} 
Applying \eqref{lambdaPPa}, \eqref{phi-psi-2} and $X_- =Y_-$ 
to one relation in \eqref{ricci3}, we obtain 
\begin{equation} 
(P_+ -Q_+ +2\Tilde{P}_- )_{st} =L_0 e^{-P_+ -\Tilde{P}_-} , 
\label{PQ2P} 
\end{equation} 
where $\Tilde{P}_- :=P_- -2a(s)$. 
Applying \eqref{PQ}, \eqref{codazzi5}, \eqref{lambdaPPa} 
to the other relation in \eqref{ricci3}, 
we obtain 
\begin{equation} 
(P_+ +Q_+ )_{st} +\varepsilon e^{P_+ +Q_+} =L_0 e^{-P_+ -\Tilde{P}_-} 
\label{P_++Q+} 
\end{equation} 
for $\varepsilon \in \{ +, -\}$. 
If we set 
\begin{equation} 
f_1 =P_+ -Q_+ +2\Tilde{P}_- , \quad 
f_2 =     Q_+ - \Tilde{P}_- , 
\label{f1f2} 
\end{equation} 
then from \eqref{PQ2P} and \eqref{P_++Q+}, we obtain 
\begin{equation} 
(f_1 )_{st} =  L_0 e^{-f_1 -f_2} , \quad 
(f_2 )_{st} =-\dfrac{\varepsilon}{2} e^{f_1 +2f_2} . 
\label{f1f2eq} 
\end{equation} 
These equations form a semilinear hyperbolic system. 

Let $f_1$, $f_2$ be functions of two variables $s$, $t$ 
satisfying \eqref{f1f2eq} for $\varepsilon \in \{ +, -\}$. 
See \cite[Chapter 5]{CH} or \cite{ando17} 
for existence and uniqueness of solutions 
of semilinear hyperbolic systems.  
Let $P_+$, $Q_+$, $\Tilde{P}_-$ be functions satisfying \eqref{f1f2}. 
Then $P_+$, $Q_+$, $\Tilde{P}_-$ satisfy \eqref{PQ2P} and \eqref{P_++Q+}. 
We set 
\begin{equation} 
\begin{split} 
& \lambda :=-\dfrac{1}{2} (P_+ +\Tilde{P}_- ), \\ 
& \phi_-  :=-\dfrac{1}{\sqrt{2}} (P_+ -Q_+ )_s -(\Tilde{P}_- )_u , \quad 
  \psi_-  :=-\dfrac{1}{\sqrt{2}} (P_+ -Q_+ )_s -(\Tilde{P}_- )_v . 
\end{split} 
\label{lpp} 
\end{equation} 
Let $\phi_+$, $\psi_+$ be as in \eqref{codazzi5}. 
Then $\lambda$, $\phi_{\pm}$, $\psi_{\pm}$ satisfy \eqref{mu0}. 
Therefore we define $\mu_1$, $\mu_2$ by 
\begin{equation} 
\begin{split} 
\mu_1 & := (\lambda_v -\psi_+ =) \ 
           -\dfrac{1}{2} (P_+ )_u 
           +\dfrac{1}{\sqrt{2}} (Q_+ )_s 
           -\dfrac{1}{2} (\Tilde{P}_- )_v , \\ 
\mu_2 & := (\lambda_u -\phi_+ =) \ 
           -\dfrac{1}{2} (P_+ )_v 
           +\dfrac{1}{\sqrt{2}} (Q_+ )_s 
           -\dfrac{1}{2} (\Tilde{P}_- )_u . 
\end{split} 
\label{mu1mu22} 
\end{equation} 
We define $X_+$, $Y_+$ by 
\begin{equation} 
X_+ 
 :=\dfrac{\varepsilon'}{2} (e^{P_+} +\varepsilon e^{Q_+} ), 
   \quad 
Y_+ 
 :=\dfrac{\varepsilon'}{2} (e^{P_+} -\varepsilon e^{Q_+} ) 
\label{setXpmYpm+} 
\end{equation} 
for $\varepsilon' \in \{ +, -\}$. 
Then we have $X^2_+ \not= Y^2_+$.  
and set $X_- =Y_- :=(\varepsilon'' /2)e^{\Tilde{P}_-}$ 
for $\varepsilon'' \in \{ +, -\}$.  
We define $\alpha_k$, $\beta_k$ ($k=1, 2$) by \eqref{alphabeta}. 
Then $S$, $T$ with \eqref{ST} satisfy $S_v -T_u =ST-TS$. 

Hence we obtain 

\begin{thm}\label{thm:H=0Qzn} 
Let $N$ be a 4-dimensional neutral space form 
with constant sectional curvature $L_0$. 
Let $M$ be a Lorentz surface 
and $F:M\rightarrow N$ a time-like and conformal immersion of $M$ into $N$ 
with zero mean curvature vector.  
Then the following are equivalent\/$:$ 
\begin{itemize}
\item[{\rm (a)}]{$F$ satisfies $X^2_+ \not= Y^2_+$ and $X_- =Y_- \not= 0;$} 
\item[{\rm (b)}]{the induced metric $g$ by $F$ is locally represented 
as $g=e^{2\lambda} d\check{w}d\overline{\check{w}}$ 
for a real-valued function $\lambda$ given by \eqref{lpp} 
and $\alpha_k$, $\beta_k$, $\mu_k$ $(k=1, 2)$ 
satisfy \eqref{alphabeta}, \eqref{mu1mu22} 
with \eqref{f1f2}, \eqref{f1f2eq}, \eqref{setXpmYpm+} 
and $X_- =Y_- :=(\varepsilon'' /2)e^{\Tilde{P}_-}$.} 
\end{itemize} 
In addition, functions $\lambda$, $\alpha_k$, $\beta_k$, $\mu_k$ $(k=1, 2)$ 
as in {\rm (b)} locally define a time-like surface in $N$ 
with zero mean curvature vector, $X^2_+ \not= Y^2_+$ and $X_- =Y_- \not= 0$, 
which is unique up to an isometry of $N$. 
\end{thm} 

\appendix 

\section{Time-like surfaces with zero mean curvature vector 
and \mbox{\boldmath{$K\equiv L_0$}} 
in a 4-dimensional Lorentzian space form 
with constant sectional curvature \mbox{\boldmath{$L_0$}}}\label{sect:Lsf} 

\setcounter{equation}{0} 

In this appendix, 
let $N$ be a 4-dimensional Lorentzian space form 
with constant sectional curvature $L_0$. 
If $L_0 =0$, 
then we can suppose $N=E^4_1$; 
if $L_0 >0$, 
then we can suppose $N=\{ x\in E^5_1 \ | \ \langle x, x\rangle =1/L_0 \}$; 
if $L_0 <0$, 
then we can suppose $N=\{ x\in E^5_2 \ | \ \langle x, x\rangle =1/L_0 \}$, 
where $\langle \ , \ \rangle$ is the metric of $E^5_1$ or $E^5_2$. 
Let $h$ be the Lorentzian metric of $N$. 
Let $M$ be a Lorentz surface 
and $F:M\rightarrow N$ a time-like and conformal immersion of $M$ into $N$ 
with zero mean curvature vector. 
Let $\check{w} =u+jv$ be a local paracomplex coordinate of $M$. 
We represent the induced metric $g$ on $M$ by $F$ 
as $g=e^{2\lambda} d\check{w}d\overline{\check{w}} 
     =e^{2\lambda} (du^2 -dv^2 )$ 
for a real-valued function $\lambda$. 
Let $\Tilde{\nabla}$ denote the connection of $E^4_1$, $E^5_1$ or $E^5_2$ 
according to $L_0 =0$, $>0$ or $<0$. 
We set $T_1 :=dF(\partial /\partial u )$, $T_2 :=dF(\partial /\partial v )$. 
Let $N_1$, $N_2$ be normal vector fields of $F$ satisfying 
$$h(N_1 , N_1 )=h(N_2 , N_2 )=e^{2\lambda} , \quad 
  h(N_1 , N_2 )=0.$$ 
Then we have 
\begin{equation} 
\begin{split} 
\Tilde{\nabla}_{T_1} (T_1 \ T_2 \ N_1 \ N_2 \ F) 
& =                  (T_1 \ T_2 \ N_1 \ N_2 \ F)S, \\ 
\Tilde{\nabla}_{T_2} (T_1 \ T_2 \ N_1 \ N_2 \ F) 
& =                  (T_1 \ T_2 \ N_1 \ N_2 \ F)T, 
\end{split} 
\label{dt1t2L} 
\end{equation} 
where 
\begin{equation} 
S =\left[ \begin{array}{ccccc} 
       \lambda_u         & \lambda_v & -\alpha_1  & -\beta_1   & 1 \\ 
       \lambda_v         & \lambda_u &  \alpha_2  &  \beta_2   & 0 \\ 
       \alpha_1          & \alpha_2  &  \lambda_u & -\mu_1     & 0 \\ 
       \beta_1           & \beta_2   &  \mu_1     &  \lambda_u & 0 \\ 
       -L_0 e^{2\lambda} &  0        &   0        &   0        & 0 
            \end{array} 
   \right] , \quad 
T =\left[ \begin{array}{ccccc} 
       \lambda_v & \lambda_u         & -\alpha_2  & -\beta_2   & 0 \\ 
       \lambda_u & \lambda_v         &  \alpha_1  &  \beta_1   & 1 \\ 
       \alpha_2  & \alpha_1          &  \lambda_v & -\mu_2     & 0 \\ 
       \beta_2   & \beta_1           &  \mu_2     &  \lambda_v & 0 \\ 
        0        &  L_0 e^{2\lambda} &   0        &   0        & 0 
            \end{array} 
  \right] , 
\label{STL} 
\end{equation} 
and $\alpha_k$, $\beta_k$, $\mu_k$ ($k=1, 2$) are real-valued functions. 
From \eqref{dt1t2L}, we obtain $S_v -T_u =ST-TS$. 
Suppose that the curvature $K$ of $g$ is identically equal to $L_0$. 
Then the equation of Gauss is given by 
\begin{equation} 
\alpha^2_1 +\beta^2_1 =\alpha^2_2 +\beta^2_2 . 
\label{gaussL}
\end{equation}
The equations of Codazzi are given by 
\begin{equation} 
\begin{split}
  (\alpha_1 )_v -(\alpha_2 )_u 
& =\alpha_2 \lambda_u -\alpha_1 \lambda_v 
  -\beta_2  \mu_1     +\beta_1   \mu_2  , \\ 
  (\alpha_2 )_v -(\alpha_1 )_u 
& =\alpha_1 \lambda_u -\alpha_2 \lambda_v 
  -\beta_1  \mu_1     +\beta_2   \mu_2  , \\ 
  (\beta_1 )_v -(\beta_2 )_u 
& =\beta_2 \lambda_u  -\beta_1 \lambda_v 
  +\alpha_2  \mu_1    -\alpha_1  \mu_2  , \\ 
  (\beta_2 )_v -(\beta_1 )_u 
& =\beta_1 \lambda_u  -\beta_2 \lambda_v 
  +\alpha_1  \mu_1    -\alpha_2  \mu_2 .  
\end{split} 
\label{codazziL} 
\end{equation} 
The equation of Ricci is given by 
\begin{equation} 
  (\mu_1 )_v -(\mu_2 )_u =2\alpha_1 \beta_2 -2\alpha_2 \beta_1 . 
\label{ricciL}
\end{equation}
We see that \eqref{ricciL} coincides with \eqref{ricci}. 
We set 
\begin{equation*} 
p_{\pm} :=e^{\lambda} (\alpha_1 \pm \alpha_2 ), \quad 
q_{\pm} :=e^{\lambda} (\beta_1  \pm \beta_2  ), \quad 
\mu_{\pm} :=\mu_1 \pm \mu_2 . 
\end{equation*} 
Then the system \eqref{codazziL} is rewritten into 
a system which consists of 
\begin{equation} 
(p_+ )_t = \dfrac{1}{\sqrt{2}} \mu_- q_+ , \quad 
(q_+ )_t =-\dfrac{1}{\sqrt{2}} \mu_- p_+ 
\label{codazziL21} 
\end{equation} 
and 
\begin{equation} 
(p_- )_s = \dfrac{1}{\sqrt{2}} \mu_+ q_- , \quad 
(q_- )_s =-\dfrac{1}{\sqrt{2}} \mu_+ p_- . 
\label{codazziL22} 
\end{equation} 
By \eqref{codazziL21}, 
we can find a function $\Tilde{\mu}_-$ satisfying $(\Tilde{\mu}_- )_t =\mu_-$ 
and that $p_+$, $q_+$ are represented as 
\begin{equation} 
p_+ = C_+ (s) \cos \left( \dfrac{\Tilde{\mu}_-}{\sqrt{2}} \right) , \quad 
q_+ =-C_+ (s) \sin \left( \dfrac{\Tilde{\mu}_-}{\sqrt{2}} \right) 
\label{codazziL31} 
\end{equation} 
for a function $C_+$ of one variable; 
by \eqref{codazziL22}, 
we can find a function $\Tilde{\mu}_+$ 
satisfying $(\Tilde{\mu}_+ )_s =\mu_+$ 
and that $p_-$, $q_-$ are represented as 
\begin{equation} 
p_- = C_- (t) \cos \left( \dfrac{\Tilde{\mu}_+}{\sqrt{2}} \right) , \quad 
q_- =-C_- (t) \sin \left( \dfrac{\Tilde{\mu}_+}{\sqrt{2}} \right) 
\label{codazziL32} 
\end{equation} 
for a function $C_-$ of one variable. 
Since we have $p_+ p_- +q_+ q_- =0$ from \eqref{gaussL}, 
we see from \eqref{codazziL31} and \eqref{codazziL32} that 
$F$ satisfies $C_+ =0$, $C_- =0$ 
or $\cos ((\Tilde{\mu}_+ -\Tilde{\mu}_- )/\sqrt{2} )=0$. 
Suppose $\cos ((\Tilde{\mu}_+ -\Tilde{\mu}_- )/\sqrt{2} )=0$. 
Then $\Tilde{\mu}_- =\Tilde{\mu}_+ +(n+1/2)\sqrt{2} \pi$ 
for an integer $n$. 
Then we have $(\Tilde{\mu}_+ )_s =(\Tilde{\mu}_- )_s =\mu_+$ 
and          $(\Tilde{\mu}_+ )_t =(\Tilde{\mu}_- )_t =\mu_-$. 
These imply $(\mu_1 )_v =(\mu_2 )_u$. 
We have 
\begin{equation*} 
 \alpha_1 \beta_2 -\alpha_2 \beta_1 
=\dfrac{1}{2e^{2\lambda}} (p_- q_+ -p_+ q_- ) 
=\dfrac{C_+ (s)C_- (t)}{2e^{2\lambda}} 
 \sin \left( \dfrac{1}{\sqrt{2}} (\Tilde{\mu}_+ -\Tilde{\mu}_- )\right) . 
\label{a1b2a2b1} 
\end{equation*} 
Since $\Tilde{\mu}_- =\Tilde{\mu}_+ +(n+1/2)\sqrt{2} \pi$, 
we obtain $\sin ((1/\sqrt{2} )(\Tilde{\mu}_+ -\Tilde{\mu}_- ))\not= 0$. 
Therefore from \eqref{ricciL} and $(\mu_1 )_v =(\mu_2 )_u$, 
we obtain $C_+ =0$ or $C_- =0$. 
Then 
\begin{equation} 
(\alpha_2 , \beta_2 )=\varepsilon (\alpha_1 , \beta_1 ) 
\label{alpha2beta2} 
\end{equation} 
for $\varepsilon \in \{ +, -\}$. 
By $(\mu_1 )_v =(\mu_2 )_u$, 
there exists a function $\gamma$ 
satisfying $\gamma_u =\mu_1$, $\gamma_v =\mu_2$. 
Therefore we have $\gamma_s =\mu_+ /\sqrt{2}$, 
                  $\gamma_t =\mu_- /\sqrt{2}$. 
Noticing $(\Tilde{\mu}_+ )_s =(\Tilde{\mu}_- )_s =\mu_+$ 
and      $(\Tilde{\mu}_+ )_t =(\Tilde{\mu}_- )_t =\mu_-$, 
we suppose     $\Tilde{\mu}_- =\sqrt{2} \gamma$ 
(respectively, $\Tilde{\mu}_+ =\sqrt{2} \gamma$) 
if $\varepsilon =+$ (respectively, $-$). 
Then we obtain 
\begin{equation} 
\begin{split} 
\alpha_1 = \dfrac{C_{\varepsilon} (u+\varepsilon v)}{2e^{\lambda}} 
           \cos \gamma , \quad 
\beta_1  =-\dfrac{C_{\varepsilon} (u+\varepsilon v)}{2e^{\lambda}} 
           \sin \gamma . 
\end{split} 
\label{alpha1beta1} 
\end{equation} 

Let $\lambda$ be a function of two variables $u$, $v$ 
satisfying \eqref{K=L0}.  
Let $\gamma$ be a function of two variables $u$, $v$. 
Let $\alpha_1$, $\beta_1$ be functions given by \eqref{alpha1beta1} 
for a function $C_{\varepsilon}$ of one variable 
and $\varepsilon \in \{ +, -\}$. 
Let $\alpha_2$, $\beta_2$ be as in \eqref{alpha2beta2}. 
Set $\mu_1 =\gamma_u$, $\mu_2 =\gamma_v$. 
Then $\lambda$, $\alpha_k$, $\beta_k$, $\mu_k$ ($k=1, 2$) 
satisfy $S_v -T_u =ST-TS$ for $S$, $T$ with \eqref{STL}. 

Hence we obtain 

\begin{thm}\label{thm:H=0K=L0in4dimL} 
Let $N$ be a 4-dimensional Lorentzian space form 
with constant sectional curvature $L_0$. 
Let $M$ be a Lorentz surface 
and $F:M\rightarrow N$ a time-like and conformal immersion of $M$ into $N$ 
with zero mean curvature vector.  
Then the following are equivalent\/$:$ 
\begin{itemize}
\item[{\rm (a)}]{the curvature $K$ of the induced metric $g$ by $F$ is 
identically equal to $L_0 ;$} 
\item[{\rm (b)}]{$g$ is locally represented 
as $g=e^{2\lambda} d\check{w}d\overline{\check{w}}$ 
for a real-valued function $\lambda$ satisfying \eqref{K=L0} 
and $\alpha_k$, $\beta_k$, $\mu_k$ $(k=1, 2)$ satisfy 
\eqref{alpha2beta2}, \eqref{alpha1beta1}, 
$\mu_1 =\gamma_u$, $\mu_2 =\gamma_v$ 
for a function $\gamma$ of two variables $u$, $v$.} 
\end{itemize} 
In addition, functions $\lambda$, 
$\alpha_k$, $\beta_k$, $\mu_k$ $(k=1, 2)$ as in {\rm (b)} 
locally define a time-like surface in $N$ 
with zero mean curvature vector and $K\equiv L_0$, 
which is unique up to an isometry of $N$. 
\end{thm}

\vspace{4mm} 

\par\noindent 
\footnotesize{Faculty of Advanced Science and Technology, 
              Kumamoto University \\ 
              2--39--1 Kurokami, Kumamoto 860--8555 Japan} 

\par\noindent  
\footnotesize{E-mail address: andonaoya@kumamoto-u.ac.jp} 

\end{document}